
\magnification=1200
\pageno=1
\baselineskip=18pt
\parindent=0pt

\font\bigbf=cmbx10 scaled\magstep1
\raggedbottom

\rm
\centerline{\bigbf Some integrals of hypergeometric functions}
 \vskip10pt
\centerline{by Andr\'as BIR\'O\footnote{}{Research partially supported by the NKFIH (National Research, Development and Innovation Office) Grants No. K104183, K109789, K119528, $\hbox{\rm ERC}_{
-}\hbox{\rm $\hbox{\rm HU}_{-}15\;118946$}$, and ERC-AdG Grant no. 321104}}
 \footnote{}{} \footnote{}{2000
}\hfill\break

\footnote{}{} \footnote{}{Key words and phrases: hypergeometric functions, integral formulas}\hfill\break

\centerline{A. R\'enyi Institute of Mathematics, Hungarian Academyof Sciences}

\centerline{ 1053 Budapest, Re\'altanoda u. 13-15., Hungary;
e-mail: biroand@renyi.hu}
 \vskip20pt

\noindent {\bf Abstract}.  We consider a certain definite integral involving the product of two classical hypergeometric functions having complicated arguments. We show in this paper the surprising fact that this integral does not depend on the parameters of the hypergeometric functions.

\noindent{\bf 1. Statement of the results}
\medskip

{\bf 1.1}. Our aim in this paper is to prove Theorem 1.1 stated below.

We write $ $$F\left(\matrix{\alpha ,\beta\cr \gamma\cr} ;z\right)$
for the Gauss hypergeometric function (instead of the notation
$F\left(\alpha ,\beta ,\gamma ;z\right)$ used in [3]).

{\bf THEOREM 1.1.} {\it Let} $0<T<S<1${\it , and let} $t$ {\it be any complex }
{\it number. Then}
$$A_t\left(S,T\right):=\int_T^S{{F\left(\matrix{2it,-2it\cr
{1\over 2}\cr}
;{{\left(1+\sqrt z\right)\left(\sqrt z-\sqrt T\right)}\over {2\left
(1-\sqrt T\right)\sqrt z}}\right)F\left(\matrix{it,-it\cr
{1\over 2}\cr}
;-{{\left(S-z\right)\left(1-z\right)}\over {\left(1-\sqrt S\right
)^2z}}\right)}\over {\left(1-z\right)\sqrt {z-T}\sqrt {S-z}}}dz\eqno
(1.1)$$
{\it equals}
$${{\pi}\over {\sqrt {1-T}\sqrt {1-S}}}.$$

{\bf 1.2.} We note that the numerator in (1.1) is symmetric in $S$
and $T$. Indeed, by applying the quadratic transformation
$$F\left(\matrix{it,-it\cr
{1\over 2}\cr}
;4w\left(1-w\right)\right)=F\left(\matrix{2it,-2it\cr
{1\over 2}\cr}
;w\right)\hbox{\rm \ for $-\infty <w\le{1\over 2}$}\eqno (1.2)$$
(this follows from [2], p 999, 9.133) for the second hypergeometric
function one has
$$F\left(\matrix{it,-it\cr
{1\over 2}\cr}
;-{{\left(S-z\right)\left(1-z\right)}\over {\left(1-\sqrt S\right
)^2z}}\right)=F\left(\matrix{2it,-2it\cr
{1\over 2}\cr}
;{{\left(1+\sqrt z\right)\left(\sqrt z-\sqrt S\right)}\over {2\left
(1-\sqrt S\right)\sqrt z}}\right).\eqno (1.3)$$
We have decided to use the left-hand side of (1.3) in
Theorem 1.1 because the argument of the hypergeometric
function there is a bit simpler than on the right-hand
side.

We also note that we could not simply apply (1.2) for
the first hypergeometric function in (1.1) because the
argument of the function there can be any number
between 0 and 1, and we cannot apply (1.2) for $w>{1\over 2}$.

{\bf 1.3.} We think that the identity of Theorem 1.1 is
interesting in its own right, but we mention that we
observed it while studying the integral operator
$$g\left(S\right)=\int_0^S{{F\left(\matrix{it,-it\cr
{1\over 2}\cr}
;-{{\left(S-z\right)\left(1-z\right)}\over {\left(1-\sqrt S\right
)^2z}}\right)}\over {\sqrt {S-z}}}f\left(z\right)dz,$$
where $f$ and $g$ are functions on $(0,1)$. In fact, one can
give the inverse of this transform using Theorem 1.1, we
intend to show it in a forthcoming paper.

We also mention that when $T$ and $S$ are fixed, but $t$ is
a real number and tends to infinity, then the integrand
in the integral defining
$A_t\left(S,T\right)$ may be exponentially large, so it is an interesting fact
that $A_t\left(S,T\right)$ itself is bounded (which follows, of
course, from Theorem 1.1).

\noindent{\bf 2. Preliminary lemmas}
\medskip

We first need a lemma which shows that $F\left(\matrix{-it,it\cr
{1\over 2}\cr}
;-x\right)$
is in fact a trigonometric function, and using this fact
we also show a product formula for this function.

For a complex number $z\neq 0$ we set its argument in $(-\pi ,\pi
]$,
and write $\log z=\log\left|z\right|+i\arg z,$ where $\log\left|z\right
|$ is real. We
define the power $z^s$ for any $s\in {\bf C}$ by $z^s=e^{s\log z}$.

{\bf LEMMA 2.1.} {\it (i) For every complex} $t$ {\it and for every real }
{\it number} $x>-1$ {\it we have that}
$$F\left(\matrix{-it,it\cr
{1\over 2}\cr}
;-x\right)={1\over 2}\left(\left(\sqrt {x+1}+\sqrt x\right)^{2it}
+\left(\sqrt {x+1}+\sqrt x\right)^{-2it}\right),\eqno (2.1)$$
{\it and one can also write it as }
$$F\left(\matrix{-it,it\cr
{1\over 2}\cr}
;-x\right)=\cos\left(2t\log\left(\sqrt {x+1}+\sqrt x\right)\right
)=\cos\left(2t\log\left(\sqrt {x+1}-\sqrt x\right)\right).\eqno (
2.2)$$
{\it (ii) For every complex} $t$ {\it and for every real numbers} $
x>0$
{\it and} $y>0$ {\it we have that the product}
$$2F\left(\matrix{-it,it\cr
{1\over 2}\cr}
;-x\right)F\left(\matrix{-it,it\cr
{1\over 2}\cr}
;-y\right)$$
{\it equals}
$$F\left(\matrix{-it,it\cr
{1\over 2}\cr}
;-\left(\sqrt x\sqrt {y+1}+\sqrt y\sqrt {x+1}\right)^2\right)+F\left
(\matrix{-it,it\cr
{1\over 2}\cr}
;-\left(\sqrt x\sqrt {y+1}-\sqrt y\sqrt {x+1}\right)^2\right).$$
{\it Proof.\/} For the case $x>0$ formula (2.1) follows from [2],
p. 998, 9.131.1 and [3], (1.5.19). Then it follows by analytic
continuation also for $x>-1$, taking into account that
$$\left(\sqrt {x+1}+\sqrt x\right)^{-2it}=\left(\sqrt {x+1}-\sqrt
x\right)^{2it}.$$
Then (2.2) follows at once, the second equality there follows by the evenness of the
cosine function.

Now, writing
$$X=\left(\sqrt x\sqrt {y+1}+\epsilon\sqrt y\sqrt {x+1}\right)^2\eqno
(2.3)$$
we have
$$X+1=\left(\sqrt {x+1}\sqrt {y+1}+\epsilon\sqrt x\sqrt y\right)^
2\eqno (2.4)$$
for $\epsilon =1$ and also for $\epsilon =-1$, hence
$$F\left(\matrix{-it,it\cr
{1\over 2}\cr}
;-X\right)=\cos\left(2t\log\left(\left(\sqrt {x+1}+\sqrt x\right)\left
(\sqrt {y+1}+\epsilon\sqrt y\right)\right)\right)$$
Part (ii) then follows by the trigonometric identity
[2], p 29, 1.314.3. The lemma is proved.

The following three lemmas are easy consequences of
Lemma 2.1 and the identity
$${1\over {2\pi}}\int_0^{\infty}{{\Gamma\left(a\pm is\right)\Gamma\left
(b\pm is\right)\Gamma\left(c\pm is\right)}\over {\Gamma\left(\pm
2is\right)}}ds=\Gamma\left(a+b\right)\Gamma\left(a+c\right)\Gamma\left
(b+c\right),\eqno (2.5)$$
which is valid for $a\ge 0$ and $b,c>0$. This is formula (3.6.1)
of [1] for positive numbers $a$,$b$ and $c$, but by taking $a
\rightarrow 0+0$ the same formula holds for $a=0$ and positive $b$ and
$c$.

{\bf LEMMA 2.2.} {\it For every} $A>-1$ {\it and for every} $T\ge 0$ {\it we have }
{\it that}
$${1\over {2\pi}}\int_0^{\infty}{{\Gamma\left(\pm is\right)\Gamma\left
({1\over 2}\pm is\right)\Gamma\left({1\over 2}+T\pm is\right)}\over {
\Gamma\left(\pm 2is\right)}}F\left(\matrix{\pm is\cr
{1\over 2}\cr}
;-A\right)ds\eqno (2.6)$$
{\it equals}
$$\Gamma\left({1\over 2}\right)\Gamma\left(1+T\right)\Gamma\left({
1\over 2}+T\right)\left(1+A\right)^{-{1\over 2}-T}.\eqno (2.7)$$
{\it Proof.\/} Using Lemma 2.1 (i) and the Stirling formula to give an upper bound it is
trivial by analytic continuation that it is enough to prove the statement for
$0<A<1$.

Using (2.5) for $0<A<1$ we have that (2.6) equals
$$\sum_{k=0}^{\infty}{{\left(-A\right)^k}\over {k!\left({1\over 2}\right
)_k}}\Gamma\left({1\over 2}+k\right)\Gamma\left(1+T\right)\Gamma\left
({1\over 2}+k+T\right),$$
and this equals (2.7) by the binomial theorem. The lemma
is proved.

{\bf LEMMA 2.3.} {\it For every} $A>-1$ {\it and for every} $r>
0$ {\it we have }
{\it that}
$${1\over {2\pi}}\int_0^{\infty}{{\Gamma\left(\pm is\right)\Gamma^
2\left({1\over 2}\pm is\right)F\left(\matrix{{1\over 2}\pm is\cr
{1\over 2}\cr}
;-r\right)}\over {\Gamma\left(\pm 2is\right)}}F\left(\matrix{\pm
is\cr
{1\over 2}\cr}
;-A\right)ds\eqno (2.8)$$
{\it equals}
$$\pi\left(1+A\right)^{{1\over 2}}{1\over {1+r+A}}.$$
{\it Proof.\/} Note first that
$$F\left(\matrix{{1\over 2}\pm is\cr
{1\over 2}\cr}
;-r\right)=\left(1+r\right)^{-{1\over 2}}\cos\left(2s\log\left(\sqrt {
r+1}+\sqrt r\right)\right)\eqno (2.9)$$
by the third line of [2], p 999, 9.131.1 and by (2.2)
above. By this formula and by (2.2) one sees by analytic
continuation that it is enough to
show the statement for fixed $A>0$ and for small enough positive $
r$. In this
case, using Lemma 2.2 we see that (2.8) equals
$$\sum_{T=0}^{\infty}{{\left(-r\right)^T}\over {T!\left({1\over 2}\right
)_T}}\Gamma\left({1\over 2}\right)\Gamma\left(1+T\right)\Gamma\left
({1\over 2}+T\right)\left(1+A\right)^{-{1\over 2}-T}.$$
By summing this geometric series the lemma follows.

{\bf LEMMA 2.4.} {\it For every} $A>-1$ {\it and for every} $r>0$ and $
B\ge 0$ {\it we have }
{\it that}
$${1\over {2\pi}}\int_0^{\infty}{{\Gamma\left(\pm is\right)\Gamma^
2\left({1\over 2}\pm is\right)F\left(\matrix{{1\over 2}\pm is\cr
{1\over 2}\cr}
;-r\right)}\over {\Gamma\left(\pm 2is\right)}}F\left(\matrix{\pm
is\cr
{1\over 2}\cr}
;-A\right)F\left(\matrix{\pm is\cr
{1\over 2}\cr}
;-B\right)ds\eqno (2.10)$$
{\it equals}
$$\pi\left(1+A\right)^{{1\over 2}}\left(1+B\right)^{{1\over 2}}{{
1+A+r+B}\over {\left(1+A+r+B\right)^2+4rAB}}.\eqno (2.11)$$
{\it The denominator in (2.11) is positive.}

{\it Proof.\/} Note first that
$$\left(1+A+r+B\right)^2+4rAB>\left(r+B\right)^2-4rB=\left(r-B\right
)^2\ge 0$$
for $r,B\ge 0$ and $A>-1$, hence the last statement is true..

We may assume $B>0$, since the $B=0$ case is proved in
Lemma 2.3.

By the third line of [2], p 999, 9.131.1 and by (ii)
of Lemma 2.1 we have that
$$2\left(1+r\right)^{1/2}F\left(\matrix{{1\over 2}\pm is\cr
{1\over 2}\cr}
;-r\right)F\left(\matrix{\pm is\cr
{1\over 2}\cr}
;-B\right)$$
equals
$$\sum_{\epsilon\in \{-1,1\}}F\left(\matrix{-is,is\cr
{1\over 2}\cr}
;-\left(\sqrt r\sqrt {B+1}+\epsilon\sqrt B\sqrt {r+1}\right)^2\right
).$$
Applying again the third line of [2], p 999, 9.131.1 and
then using Lemma 2.3, taking into account (2.3) and (2.4)
with $r$ and $B$ in place of $x$ and $y$, we get that (2.10) equals
$${{\pi\left(1+A\right)^{{1\over 2}}}\over {2\left(1+r\right)^{1/
2}}}\sum_{\epsilon\in \{-1,1\}}{{\sqrt {r+1}\sqrt {B+1}+\epsilon\sqrt
r\sqrt B}\over {1+\left(\sqrt r\sqrt {B+1}+\epsilon\sqrt B\sqrt {
r+1}\right)^2+A}}.$$
Then by the identities
$$\left(1+A+r+B+2rB\right)^2-4r\left(1+r\right)B\left(1+B\right)=\left
(1+A+r+B\right)^2+4rAB$$
and
$$\sum_{\epsilon\in \{-1,1\}}\left(\sqrt {r+1}\sqrt {B+1}+\epsilon\sqrt
r\sqrt B\right)\left(1+\left(\sqrt r\sqrt {B+1}-\epsilon\sqrt B\sqrt {
r+1}\right)^2+A\right)=$$
$$=2\sqrt {r+1}\sqrt {B+1}\left(1+A+r+B\right)$$
we get the lemma.

\noindent{\bf 3. Proof of the theorem}
\medskip

We now turn to the proof of Theorem 1.1. Let
$0<T<S<1$ be given from now on. We note that $A_t\left(S,T\right)$ is obviously
an entire function of $t$, hence it is enough to prove the
identity for real $t$.

{\bf LEMMA 3.1.} {\it Assume that there is an} $r_0>0$ {\it such that}
$${1\over {\pi}}\int_0^{\infty}{{\Gamma\left(\pm 2it\right)\Gamma^
2\left({1\over 2}\pm 2it\right)F\left(\matrix{{1\over 2}\pm 2it\cr
{1\over 2}\cr}
;-r\right)}\over {\Gamma\left(\pm 4it\right)}}\left(A_t\left(S,T\right
)-{{\pi}\over {\sqrt {1-T}\sqrt {1-S}}}\right)dt=0\eqno (3.1)$$
{\it for every} $r>r_0${\it . Then }
$$A_t\left(S,T\right)={{\pi}\over {\sqrt {1-T}\sqrt {1-S}}}$$
{\it for every real} $t${\it .}

{\it Proof.\/} By Lemma 2.1 (i) for
$0<Y<1$ one has that
$$F\left(\matrix{2it,-2it\cr
{1\over 2}\cr}
;Y\right)=\cos\left(4it\arg\left(\sqrt {1-Y}+i\sqrt Y\right)\right
),$$
where the argumentum of a complex number lying in the
right half-plane is taken in $\left(-{{\pi}\over 2},{{\pi}\over 2}\right
).$ Estimating the first
hypergeometric function in (1.1) by this relation and
estimating the second one directly by (2.2), we get that
$$A_t\left(S,T\right)=O_{S,T}\left(e^{\left(2\pi -c\left(S,T\right
)\right)t}\right)\eqno (3.2)$$
with some constant $c\left(S,T\right)>0$ depending only on $S$ and $
T$.

Using (3.2), (2.9) and the Stirling formula we see that the left-hand side
of (3.1) is absolutely convergent for real $r$, moreover,
this integral extends as a holomorphic function of $r$ to
a domain in the complex plane containing the positive
real line. Therefore by the unicity theorem we see
that (3.1) is true for every $r>0$. Using (2.9) one then sees
that the Fourier transform of the function
$${{\Gamma\left(\pm 2it\right)\Gamma^2\left({1\over 2}\pm 2it\right
)}\over {\Gamma\left(\pm 4it\right)}}\left(A_t\left(S,T\right)-{{
\pi}\over {\sqrt {1-T}\sqrt {1-S}}}\right)$$
is identically 0. The lemma follows.

Hence for the proof of the theorem it is enough to show
(3.1) for large enough $r$. The next two lemma will be useful to
compute the left-hand side of (3.1).

{\bf LEMMA 3.2.} {\it For every} $r>0$ {\it one has that}

$${1\over {\pi}}\int_0^{\infty}{{\Gamma\left(\pm 2it\right)\Gamma^
2\left({1\over 2}\pm 2it\right)F\left(\matrix{{1\over 2}\pm 2it\cr
{1\over 2}\cr}
;-r\right)}\over {\Gamma\left(\pm 4it\right)}}{{\pi}\over {\sqrt {
1-T}\sqrt {1-S}}}dt$$
{\it equals }
$$\left(1+r\right)^{-1}{{\pi^2}\over {\sqrt {1-T}\sqrt {1-S}}}.$$
{\it Proof.\/} By the third line of [2], p 999, 9.131.1 and by
Lemma 2.2 with $T=0$ and $A=r$, using the substitution
$s=2t$ we get the statement.

{\bf LEMMA 3.3.} {\it For} $T\le z\le S$ {\it let }
$$B_t\left(z,S,T\right):=F\left(\matrix{2it,-2it\cr
{1\over 2}\cr}
;{{\left(1+\sqrt z\right)\left(\sqrt z-\sqrt T\right)}\over {2\left
(1-\sqrt T\right)\sqrt z}}\right)F\left(\matrix{it,-it\cr
{1\over 2}\cr}
;-{{\left(S-z\right)\left(1-z\right)}\over {\left(1-\sqrt S\right
)^2z}}\right).$$
{\it Then for every} $r>0$ {\it and for every} $T\le z\le S$ {\it we have that}
$${1\over {\pi}}\int_0^{\infty}{{\Gamma\left(\pm 2it\right)\Gamma^
2\left({1\over 2}\pm 2it\right)F\left(\matrix{{1\over 2}\pm 2it\cr
{1\over 2}\cr}
;-r\right)B_t\left(z,S,T\right)}\over {\Gamma\left(\pm 4it\right)}}
dt=I_1\left(z\right)I_2\left(r,z\right)\eqno (3.3)$$
{\it with the abbreviations}
$$I_1\left(z\right):=\pi\left(1-\sqrt z\right)\left(\sqrt z+\sqrt
T\right)^{{1\over 2}}\left(\sqrt z+\sqrt S\right)^{{1\over 2}}\left
(1-\sqrt T\right)^{{1\over 2}}\left(1-\sqrt S\right)^{{1\over 2}}$$
{\it and }
$$I_2\left(r,z\right):={{\sqrt T+\sqrt S-2\sqrt S\sqrt T+R\sqrt z
+z\left(\sqrt T+\sqrt S-2\right)}\over {E\left(r\right)+F\left(r\right
)z+G\left(r\right)z^2}}\eqno (3.4)$$
{\it with the notations (note that R also depends on} $r${\it , }
{\it but for simplicity we do not denote it)}
$$R:=2r\left(1-\sqrt T\right)\left(1-\sqrt S\right),$$
$$E\left(r\right):=\left(\sqrt T+\sqrt S-2\sqrt S\sqrt T\right)^2
+2R\sqrt S\sqrt T,$$
$$F\left(r\right):=2\left(\sqrt T+\sqrt S-2\right)\left(\sqrt T+\sqrt
S-2\sqrt S\sqrt T\right)+2R\left(1+\sqrt S\sqrt T-2\sqrt T-2\sqrt
S\right)+R^2,$$
$$G\left(r\right):=\left(\sqrt T+\sqrt S-2\right)^2+2R.$$
{\it For} $T\le z\le S$ {\it we have }
$$E\left(r\right)+F\left(r\right)z+G\left(r\right)z^2>0.\eqno (3.5)$$
{\it Proof.\/} By [2], p 999, 9.133 we have that
$$F\left(\matrix{it,-it\cr
{1\over 2}\cr}
;-{{\left(S-z\right)\left(1-z\right)}\over {\left(1-\sqrt S\right
)^2z}}\right)=F\left(\matrix{2it,-2it\cr
{1\over 2}\cr}
;-{{\left(1+\sqrt z\right)\left(\sqrt S-\sqrt z\right)}\over {2\left
(1-\sqrt S\right)\sqrt z}}\right)$$
by the conditions $0<T\le z\le S<1.$

Applying Lemma 2.4 with the substitution $s=2t$ and with
$$A=A\left(z\right):=-{{\left(1+\sqrt z\right)\left(\sqrt z-\sqrt
T\right)}\over {2\left(1-\sqrt T\right)\sqrt z}},\qquad B=B\left(
z\right):={{\left(1+\sqrt z\right)\left(\sqrt S-\sqrt z\right)}\over {
2\left(1-\sqrt S\right)\sqrt z}}\eqno (3.6)$$
($A>-1$ and $B\ge 0$ are satisfied by our conditions), using
the easily proved identities (recall the notations given in the text of
the lemma)
$$1+A\left(z\right)={{\left(1-\sqrt z\right)\left(\sqrt z+\sqrt T\right
)}\over {2\left(1-\sqrt T\right)\sqrt z}},\qquad 1+B\left(z\right
)={{\left(1-\sqrt z\right)\left(\sqrt z+\sqrt S\right)}\over {2\left
(1-\sqrt S\right)\sqrt z}},$$
$$1+A\left(z\right)+r+B\left(z\right)={{\sqrt T+\sqrt S-2\sqrt S\sqrt
T+R\sqrt z+z\left(\sqrt T+\sqrt S-2\right)}\over {2\left(1-\sqrt
T\right)\left(1-\sqrt S\right)\sqrt z}},\eqno (3.7)$$

$$\left(1+A\left(z\right)+r+B\left(z\right)\right)^2+4rA\left(z\right
)B\left(z\right)={{E\left(r\right)+F\left(r\right)z+G\left(r\right
)z^2}\over {4z\left(1-\sqrt T\right)^2\left(1-\sqrt S\right)^2}}\eqno
(3.8)$$
we obtain (3.3). The inequality (3.5) follows from (3.8)
and the last statement of Lemma 2.4. The lemma is
proved.

{\bf LEMMA 3.4.} {\it Assume that there is an} $r_0>0$ {\it such that}
$$Q_{S,T}\left(r\right):=\left(1+r\right)\int_T^S{{I_2\left(r,z\right
)}\over {\left(1+\sqrt z\right)\sqrt {\sqrt z-\sqrt T}\sqrt {\sqrt
S-\sqrt z}}}dz\eqno (3.9)$$
{\it equals}
$${{\pi}\over {\sqrt {1-T}\sqrt {1-S}\sqrt {1-\sqrt T}\sqrt {1-\sqrt
S}}}\eqno (3.10)$$
{\it for every} $r>r_0${\it , where} $I_2\left(r,z\right)$ {\it is defined in Lemma 3.3. Then Theorem 1.1 is true.}

{\it Proof.\/} This follows at once from Lemma 3.1, Lemma 3.2,
(1.1) and Lemma 3.3.

{\bf LEMMA 3.5.} {\it We have that}
$$Q_{S,T}\left(r\right)-{{\pi}\over {\sqrt {1-T}\sqrt {1-S}\sqrt {
1-\sqrt T}\sqrt {1-\sqrt S}}}=O_{S,T}\left({1\over r}\right)$$
{\it as} $r\rightarrow +\infty${\it .}

{\it Proof.\/} By (3.4), (3.7), (3.8) we have the identity
$$\left(1+r\right)I_2\left(r,z\right)={1\over {2\left(1-\sqrt T\right
)\left(1-\sqrt S\right)\sqrt z}}{{\left(1+r\right)\left(1+A\left(
z\right)+r+B\left(z\right)\right)}\over {\left(1+A\left(z\right)+
r+B\left(z\right)\right)^2+4rA\left(z\right)B\left(z\right)}},$$
so $\left(1+r\right)I_2\left(r,z\right)$ equals
$${1\over {2\left(1-\sqrt T\right)\left(1-\sqrt S\right)\sqrt z}}
+O_{S,T}\left({1\over r}\right)$$
as $r\rightarrow +\infty$, uniformly for $T\le z\le S$. By (3.9) we see that
for the proof of the present lemma it is enoiugh to show
that
$$\int_T^S{1\over {\left(1+\sqrt z\right)\sqrt {\sqrt z-\sqrt T}\sqrt {\sqrt
S-\sqrt z}}}{{dz}\over {2\sqrt z}}={{\pi\sqrt {1-\sqrt T}\sqrt {1
-\sqrt S}}\over {\sqrt {1-T}\sqrt {1-S}}}\eqno .(3.11)$$
By the substitution
$$q={{\sqrt z-\sqrt T}\over {\sqrt S-\sqrt T}},$$
we see that for (3.11) it is enough to show that
$$\int_0^1{1\over {\left(1+\sqrt T+q\left(\sqrt S-\sqrt T\right)\right
)\sqrt q\sqrt {1-q}}}dq={{\pi\sqrt {1-\sqrt T}\sqrt {1-\sqrt S}}\over {\sqrt {
1-T}\sqrt {1-S}}}.\eqno (3.12)$$
Now, the left-hand side of (3.12) equals
$${{\pi}\over {\left(1+\sqrt T\right)}}F\left(\matrix{{1\over 2},
1\cr
1\cr}
;-{{\sqrt S-\sqrt T}\over {1+\sqrt T}}\right)$$
by [2], p 995, 9.111, and this shows (3.12) by the
binomial theorem. The lemma follows.

{\bf LEMMA 3.6.} {\it For large enough positive} $r$ {\it we have that}
$$Q_{S,T}\left(r\right)-{{\pi}\over {\sqrt {1-T}\sqrt {1-S}\sqrt {
1-\sqrt T}\sqrt {1-\sqrt S}}}$$
{\it equals}
$$2in_r\pi\sqrt {\left(1-\sqrt T\right)\left(1-\sqrt S\right)}{{\sqrt
r\left(1+r\right)}\over {\sqrt {F\left(r\right)^2-4E\left(r\right
)G\left(r\right)}}},$$
{\it where} $-4\le n_r\le 4$ {\it is an integer depending on} $r${\it .}

{\bf REMARK.} The integer $n_r$ may depend also on $S$ and $T$,
but $S$ and $T$ are fixed, so we do not denote it.

{\it Proof.\/} Since $F\left(r\right)^2-4E\left(r\right)G\left(r\right
)$ is a polynomial in $r$ of exact
order $4$ and $E\left(r\right)+F\left(r\right)+G$$\left(r\right)$ is a polynomial in $
r$ of exact order
$2$ we see (since obviously $E\left(r\right)>0$ and $G\left(r\right
)>0$) that if
$r$ is large enough, then the polynomial
$$E\left(r\right)+F\left(r\right)z+G\left(r\right)z^2$$
(we consider it as a polynomial in $z$) is of exact order $2$, it does not have a double
root in $z$ and $z=1$ is not a root, $z=0$ is not a root. One has
$$E\left(r\right)+F\left(r\right)z+G\left(r\right)z^2=E\left(r\right
)\left(1-\alpha_1z\right)\left(1-\alpha_2z\right)\eqno (3.13)$$
with some complex numbers $\alpha_1$ and $\alpha_2$ different from
each other and from $1$ and $0$. Note that $\alpha_1$ and $\alpha_
2$ depend
on $r$, but for simplicity we do not denote it.

Then it is easy to see (since writing $y=\sqrt z$ we have in
(3.14) below the quotient of a polynomial of degree 3 and a polynomial of
degree 5 in $y$, and the denominator has no double root, so
we can take a partial quotient decomposition) that
$$\sqrt z{{\sqrt T+\sqrt S-2\sqrt S\sqrt T+R\sqrt z+z\left(\sqrt
T+\sqrt S-2\right)}\over {\left(1+\sqrt z\right)\left(E\left(r\right
)+F\left(r\right)z+G\left(r\right)z^2\right)}}\eqno (3.14)$$
equals
$${a\over {1+\sqrt z}}+{b\over {1-\sqrt {\alpha_1}\sqrt z}}+{c\over {
1+\sqrt {\alpha_1}\sqrt z}}+{d\over {1-\sqrt {\alpha_2}\sqrt z}}+{
e\over {1+\sqrt {\alpha_2}\sqrt z}}\eqno (3.15)$$
for every $y=\sqrt z$ with some complex numbers $a$, $b$, $c$, $d$ and $
e$
depending on $r$, where $\sqrt {\alpha_1}$ and $\sqrt {\alpha_2}$ are chosen arbitrarily,
but they are fixed. Note that the denominators in (3.15) are
nonzero for $T\le z\le S$ using (3.13) and (3.5).

By the substitution $q={{\sqrt z-\sqrt T}\over {\sqrt S-\sqrt T}}$ we see by (3.9), (3.4) (3.14) and
(3.15) that $Q_{S,T}\left(r\right)$ equals
$$2\left(1+r\right)\int_0^1{{\left({a\over {1+\sqrt z}}+{b\over {
1-\sqrt {\alpha_1}\sqrt z}}+{c\over {1+\sqrt {\alpha_1}\sqrt z}}+{
d\over {1-\sqrt {\alpha_2}\sqrt z}}+{e\over {1+\sqrt {\alpha_2}\sqrt
z}}\right)}\over {\sqrt q\sqrt {1-q}}}dq,$$
where (by a slight abuse of notation) $\sqrt z$ now denotes a
function of $q$, namely
$$\sqrt z:=\sqrt T+q\left(\sqrt S-\sqrt T\right).$$

Now, let $u$ be any such complex number for which
$$1-u\left(\sqrt T+q\left(\sqrt S-\sqrt T\right)\right)$$
is nonzero for $0\le q\le 1$. Then
$$\int_0^1{1\over {\sqrt q\sqrt {1-q}}}{1\over {\left(1-u\left(\sqrt
T+q\left(\sqrt S-\sqrt T\right)\right)\right)}}dq$$
equals
$${{\pi}\over {1-u\sqrt T}}F\left(\matrix{{1\over 2},1\cr
1\cr}
;u{{\sqrt S-\sqrt T}\over {1-u\sqrt T}}\right)={{\Gamma^2\left({1\over
2}\right)}\over {1-u\sqrt T}}\left({{1-u\sqrt S}\over {1-u\sqrt T}}\right
)^{-{1\over 2}},$$
where we used also the binomial theorem. Hence $Q_{S,T}\left(r\right
)$ equals the sum of
$$2\left(1+r\right)a{{\pi}\over {\left(1+\sqrt T\right)}}\left({{
1+\sqrt S}\over {1+\sqrt T}}\right)^{-{1\over 2}},\eqno (3.16)$$
$$D_1:=2\left(1+r\right)b{{\pi}\over {\left(1-\sqrt {\alpha_1}\sqrt
T\right)}}\left({{1-\sqrt {\alpha_1}\sqrt S}\over {1-\sqrt {\alpha_
1}\sqrt T}}\right)^{-{1\over 2}},\eqno (3.17)$$
$$D_2:=2\left(1+r\right)c{{\pi}\over {\left(1+\sqrt {\alpha_1}\sqrt
T\right)}}\left({{1+\sqrt {\alpha_1}\sqrt S}\over {1+\sqrt {\alpha_
1}\sqrt T}}\right)^{-{1\over 2}},\eqno (3.18)$$
$$D_3:=2\left(1+r\right)d{{\pi}\over {\left(1-\sqrt {\alpha_2}\sqrt
T\right)}}\left({{1-\sqrt {\alpha_2}\sqrt S}\over {1-\sqrt {\alpha_
2}\sqrt T}}\right)^{-{1\over 2}},\eqno (3.19)$$
$$D_4:=2\left(1+r\right)e{{\pi}\over {\left(1+\sqrt {\alpha_2}\sqrt
T\right)}}\left({{1+\sqrt {\alpha_2}\sqrt S}\over {1+\sqrt {\alpha_
2}\sqrt T}}\right)^{-{1\over 2}}.\eqno (3.20)$$
Note that using (3.8) we see that if
$$E\left(r\right)+F\left(r\right)z+G\left(r\right)z^2=0,$$
then
$$\left(1+A\left(z\right)+r+B\left(z\right)\right)^2+4rA\left(z\right
)B\left(z\right)=0,$$
so
$$A\left(z\right)B\left(z\right)=-{{\left(1+A\left(z\right)+r+B\left
(z\right)\right)^2}\over {4r}},$$
hence by (3.6) we have
$$\left(\sqrt z-\sqrt T\right)\left(\sqrt z-\sqrt S\right)=-{{z\left
(1-\sqrt T\right)\left(1-\sqrt S\right)}\over {\left(1+\sqrt z\right
)^2}}{{\left(1+A\left(z\right)+r+B\left(z\right)\right)^2}\over r}
.\eqno (3.21)$$
Let us determine now the coefficients $a,b,c,d,e$. For
this sake we note that (3.7) implies that
$$\sqrt T+\sqrt S-2\sqrt S\sqrt T+R\sqrt z+z\left(\sqrt T+\sqrt S
-2\right)$$
equals
$$2\left(1-\sqrt T\right)\left(1-\sqrt S\right)\sqrt z\left(1+A\left
(z\right)+r+B\left(z\right)\right),$$
so by the equality of (3.14) and (3.15), using also (3.13) we get that
$$a=-{{\left(\sqrt T+\sqrt S-2\sqrt S\sqrt T\right)-R+\left(\sqrt
T+\sqrt S-2\right)}\over {E\left(r\right)+F\left(r\right)+G\left(
r\right)}},$$
$$b={{\left(1-\sqrt T\right)\left(1-\sqrt S\right)\left(1+A\left({
1\over {\alpha_1}}\right)+r+B\left({1\over {\alpha_1}}\right)\right
)}\over {\left(\alpha_1-\alpha_2\right)\left(1+1/\sqrt {\alpha_1}\right
)E\left(r\right)}},\eqno (3.22)$$
$$c={{\left(1-\sqrt T\right)\left(1-\sqrt S\right)\left(1+A\left({
1\over {\alpha_1}}\right)+r+B\left({1\over {\alpha_1}}\right)\right
)}\over {\left(\alpha_1-\alpha_2\right)\left(1-1/\sqrt {\alpha_1}\right
)E\left(r\right)}},\eqno (3.23)$$
$$d={{\left(1-\sqrt T\right)\left(1-\sqrt S\right)\left(1+A\left({
1\over {\alpha_2}}\right)+r+B\left({1\over {\alpha_2}}\right)\right
)}\over {\left(\alpha_2-\alpha_1\right)\left(1+1/\sqrt {\alpha_2}\right
)E\left(r\right)}},\eqno (3.24)$$
$$e={{\left(1-\sqrt T\right)\left(1-\sqrt S\right)\left(1+A\left({
1\over {\alpha_2}}\right)+r+B\left({1\over {\alpha_2}}\right)\right
)}\over {\left(\alpha_2-\alpha_1\right)\left(1-1/\sqrt {\alpha_2}\right
)E\left(r\right)}}.\eqno (3.25)$$
By the definition of $R$, $E\left(r\right)$, $F\left(r\right)$ and $
G\left(r\right)$ in Lemma 3.3 we easily see that
$$a={1\over {2\left(1-\sqrt T\right)\left(1-\sqrt S\right)\left(1
+r\right)}},$$
hence (3.16) equals (3.10). Consequently
$$Q_{S,T}\left(r\right)-{{\pi}\over {\sqrt {1-T}\sqrt {1-S}\sqrt {
1-\sqrt T}\sqrt {1-\sqrt S}}}=D_1+D_2+D_3+D_4,\eqno (3.26)$$
see (3.17), (3.18),  (3.19) and  (3.20).

We claim that
$$D_i^2=-4\pi^2\left(1-\sqrt T\right)\left(1-\sqrt S\right){{r\left
(1+r\right)^2}\over {\left(\alpha_1-\alpha_2\right)^2E^2\left(r\right
)}}\eqno (3.27)$$
for every $1\le i\le 4$. For $i=1$ this follows from (3.22) and
from (3.21) with $\sqrt z={1\over {\sqrt {\alpha_1}}}$, for $i=2$ this follows from (3.23)
and from (3.21) with $\sqrt z=-{1\over {\sqrt {\alpha_1}}}$, for $
i=3$ this follows from
(3.24) and from (3.21) with $\sqrt z={1\over {\sqrt {\alpha_2}}}$, for $
i=4$ this follows from
(3.25) and from (3.21) with $\sqrt z=-{1\over {\sqrt {\alpha_2}}}$.

We get from (3.27) for every $1\le i\le 4$ that
$$D_i=2i\epsilon_i\pi\sqrt {\left(1-\sqrt T\right)\left(1-\sqrt S\right
)}{{\sqrt r\left(1+r\right)}\over {\left(\alpha_1-\alpha_2\right)
E\left(r\right)}},\eqno (3.28)$$
where $\epsilon_i=1$ or $\epsilon_i=-1$ for every $i$. Since ${1\over {
\alpha_1}}$ and ${1\over {\alpha_2}}$ are
the two roots of $E\left(r\right)+F\left(r\right)z+G\left(r\right
)z^2=0$, so
$$\alpha_1\alpha_2={{G\left(r\right)}\over {E\left(r\right)}}$$
and
$${1\over {\alpha_2}}-{1\over {\alpha_1}}=\delta{{\sqrt {F^2\left
(r\right)-4E\left(r\right)G\left(r\right)}}\over {G\left(r\right)}}
,$$
where $\delta =1$ or $\delta =-1$. Then we get
$$\left(\alpha_1-\alpha_2\right)E\left(r\right)=\delta\sqrt {F^2\left
(r\right)-4E\left(r\right)G\left(r\right)}.\eqno (3.29)$$
By (3.26), (3.28) and (3.29) we get the lemma.

{\bf LEMMA 3.7.} {\it If} $r$ {\it is a large enough positive number, then}
{\it we have that}
$$Q_{S,T}\left(r\right)={{\pi}\over {\sqrt {1-T}\sqrt {1-S}\sqrt {
1-\sqrt T}\sqrt {1-\sqrt S}}}.$$
{\it Proof.\/} It is enough to show that the integer $n_r$ in
Lemma 3.6 is $0$ for large enough $r$.

Since  $F^2\left(r\right)-4E\left(r\right)G\left(r\right)$ is a polynomial in $
r$ of exact
order $4$, therefore Lemma 3.6 implies that
$$\left|Q_{S,T}\left(r\right)-{{\pi}\over {\sqrt {1-T}\sqrt {1-S}\sqrt {
1-\sqrt T}\sqrt {1-\sqrt S}}}\right|\gg_{S,T}{{\left|n_r\right|}\over {\sqrt
r}}$$
for large $r$. By Lemma 3.5 this implies for large
enough $r$ that $n_r=0$. The lemma is proved.

Using Lemmas 3.4 and 3.7 we see that the proof of the
theorem is now complete.

 \vskip10pt

\bigskip\noindent {\bf References}

\nobreak
\parindent=12pt
\nobreak

\item{[1]} G.E. Andrews, R. Askey, R. Roy, {Special
Functions}, {\it Cambridge Univ. Press,\/} 1999

\item{[2]} I.S. Gradshteyn, I.M. Ryzhik, {Table of integrals, series and
products, 6th edition,} {\it Academic Press}, 2000

\item{[3]} L.J. Slater, {Generalized hypergeometric functions,} {\it Cambridge Univ. Press},
1966

\bye